\documentclass{ifacconf}

\usepackage{graphicx}      
\usepackage{natbib}        

\usepackage[utf8]{inputenc}

\usepackage{graphicx}
\usepackage{amsmath,bm}
\usepackage[version=4]{mhchem}
\usepackage{siunitx}
\usepackage{caption}
\usepackage{longtable,tabularx}
\usepackage{mathrsfs,amsfonts,graphicx,epsfig,stmaryrd}

\usepackage{subcaption}
\usepackage{caption}
\usepackage{color,multirow,rotating}
\usepackage{bm}
\usepackage{algorithm,algpseudocode,algorithmicx}
\usepackage{cite,url,framed,bm,balance,dsfont,varwidth}
\usepackage{nicefrac}

\newtheorem{theorem}{Theorem}

\newtheorem{proposition}{Proposition}

\newcommand{\differential}{{\rm{d}}}
\begin{document}
\begin{frontmatter}

\title{Solution of the Probabilistic Lambert's Problem: Optimal Transport Approach\thanksref{footnoteinfo}} 

\thanks[footnoteinfo]{This research was supported in part by NSF grant 2112755.}

\author[First]{Alexis M.H. Teter} 
\author[Second]{Iman Nodozi} 
\author[Third]{Abhishek Halder}

\address[First]{Department of Applied Mathematics, University of
California, Santa Cruz, CA 95064, USA (e-mail: amteter@ucsc.edu).}
\address[Second]{Department of Electrical and Computer Engineering, University of
California, Santa Cruz, CA 95064, USA (e-mail: inodozi@ucsc.edu).}
\address[Third]{Department of Aerospace Engineering, Iowa State University, IA 50011, USA (email: ahalder@iastate.edu)}

\begin{abstract}                
The deterministic variant of the Lambert's problem was posed by Lambert in the 18th century and its solution for conic trajectory has been derived by many, including Euler, Lambert, Lagrange, Laplace, Gauss and Legendre. The solution amounts to designing velocity control for steering a spacecraft from a given initial to a given terminal position subject to gravitational potential and flight time constraints. In recent years, a probabilistic variant of the Lambert's problem has received attention in the aerospace community where the endpoint position constraints are softened to endpoint joint probability distributions over the respective positions. Such probabilistic specifications account for the estimation errors, modeling uncertainties, etc. Building on a deterministic optimal control reformulation via analytical mechanics, we show that the probabilistic Lambert's problem is a generalized dynamic optimal mass transport problem where the gravitational potential plays the role of an additive state cost. This allows us to rigorously prove the existence-uniqueness of the solution for the probabilistic Lambert problem both with and without process noise. In the latter case, the problem and its solution correspond to a generalized Schrödinger bridge, much like how classical Schrodinger bridge can be seen as stochastic regularization of the optimal mass transport. We deduce the large deviation principle enjoyed by the Lambertian Schrödinger bridge. Leveraging these newfound connections, we design a computational algorithm to illustrate the nonparametric numerical solution of the probabilistic Lambert's problem.
\end{abstract}

\begin{keyword}
Lambert's Problem, Optimal Mass Transport, Schr\"{o}dinger Bridge.
\end{keyword}

\end{frontmatter}

\section{Introduction}
\textbf{Lambert's problem.} The Lambert's problem involves computing a velocity field, $\bm{v}=\bm{v}(\bm{r},t)$ subject to a two-body gravitational potential force field, that satisfies given endpoint position and hard flight time constraints:
\begin{align}
&\ddot{\bm{r}} = -\nabla_{\bm{r}}V\left(\bm{r}\right),\nonumber\\
&\bm{r}(t=t_0) = \bm{r}_0\:\text{(given)}, \quad \bm{r}(t=t_1) = \bm{r}_{1}\:\text{(given)},    \label{TPBVP}    \end{align}
where the position vector $\bm{r}=(x, y, z)^{\top}\in \mathbb{R}^3$ and the velocity vector $\bm{v}:=\dot{\bm{r}}\in \mathbb{R}^3$ are measured w.r.t. the Earth Centered Inertial (ECI) frame. The times $t_0$ and $t_1$ are the fixed initial and terminal times. The symbol $\nabla_{\bm{r}}$ denotes the Euclidean gradient w.r.t. vector $\bm{r}$. 

One can think of \eqref{TPBVP} as a \emph{partially specified} two-point boundary value problem since only the position $\bm{r}$, but not the velocity $\dot{\bm{r}}$, are prescribed at $t_0$ and $t_1$.

Let $\vert \bm{r}\vert := \sqrt{x^2 + y^2 + z^2}$. In its domain of definition, the nonlinear potential $V(\cdot)$ in \eqref{TPBVP} is assumed to be $C^{1}(\mathbb{R}^3)$, $V\in(-\infty,0)$, and is usually taken as
\begin{align}
V(\bm{r}) := - \frac{\mu}{|\bm{r}|} - \frac{\mu J_2 R_{\textrm{Earth}}^2}{2 |\bm{r}|^3} \left( 1 - \frac{3z^2}{|\bm{r}|^2} \right), 
\label{defPotential}    
\end{align}
where $\bm{r}\in [R_{\rm{Earth}}, +\infty)$, $\mu=398600.4415\ \text{km}^3/\text{s}^2$ represents the product of the Earth's gravitational constant and mass, $J_2 = 1.08263 \times 10^{-3}$ is the unitless second zonal harmonic coefficient, and $R_{\textrm{Earth}}=6378.1363$ km is the Earth's radius. When the $J_2$ term in \eqref{defPotential} is dropped, the Lambert's problem reduces to a Keplerian orbit transfer problem; see e.g., \citep[Ch. 7]{battin1999introduction}. The arc connecting $\bm{r}_0$ and $\bm{r}_1$, in the Keplerian case, is called a \emph{Keplerian arc}. In general, determining the connecting arc comprises a solution of the Lambert's problem.

In this work, we will not consider multi-revolution arcs, as was done historically for the classical Lambert's problem.

\textbf{Brief history.} The mathematics of the Lambert's problem has received interests since the dawn of analytical mechanics and astrodynamics. Historically, early motivation was determining the orbit of a comet from observations. The problem bears the name of Johann Heinrich Lambert (1728--1777), who is credited \citep[p. 141]{berggren1997pi} for the first proof of the irrationality of $\pi$ in 1761. In a letter sent to Euler in the February of the same year \citep{euler1924briefwechsel}, Lambert mentions of his result relating a given flight time $t_1-t_0$ to the problem data $\bm{r}_0,\bm{r}_1$ assuming the Keplerian arc is a parabola. Interestingly, the solutions for both the parabolic and elliptic arcs were already obtained in 1743 by Euler \citep{euler1743determinatio}. 

In his book \citep{lambert1761ih} published also in 1761, Lambert gives derivations for the parabolic, elliptic and hyberbolic cases, mentions Euler's book \citep{euler1744theoria} but not the earlier article \citep{euler1743determinatio}. In the following month, Lambert sent his book to Euler and received high praises for his results in three response letters \citep{euler1924briefwechsel}. His results in \citep{lambert1761ih} also received positive comments from Lagrange who derived several alternative proofs for such results in 1780 \citep{lagrange1780maniere}. For a detailed chronology of related results, including those by Laplace, Gauss and Legendre, we refer the readers to \citep[Sec. 9]{albouy2019lambert}. 

Lambert's problem remains relevant in modern times due to its application in trajectory design for interplanetary transfer, missile interception and rendezvous problems.

\textbf{Probabilistic Lambert's problem.} In this work, we consider the probabilistic Lambert's problem, which is a feasibility problem of the form
\begin{subequations}
\begin{align}
&\underset{\dot{\bm{r}}=\bm{v}(\bm{r},t)}{\text{find}}\quad \bm{v}\\
&\ddot{\bm{r}} = -\nabla_{\bm{r}}V\left(\bm{r}\right),\label{ProbLambertDetDyn}\\
& \bm{r}(t=t_0) \sim \rho_0\:\text{(given)}, \quad \bm{r}(t=t_1) \sim \rho_1\:\text{(given)},\label{ProbLambertEndpoint}
\end{align}
\label{ProbLambertFeasibility}
\end{subequations}
where the endpoint relative positions are random vectors with known joint probability density functions (PDFs) $\rho_0,\rho_1$. Thus, instead of navigating strictly between two specified position vectors, here the idea is to steer between their specified statistical laws.  

We assume that $\rho_0,\rho_1\in\mathcal{P}_2\left(\mathbb{R}^{3}\right)$, where $\mathcal{P}_2\left(\mathbb{R}^3\right)$ denotes set of probability density functions supported over $\mathbb{R}^{3}$ with finite second moments.

The motivation behind problem \eqref{ProbLambertFeasibility} is to allow for stochastic uncertainties in both start and endpoint position vectors. For example, the PDF $\rho_0$ may encode uncertainties due to statistical estimation errors. The PDF $\rho_1$ may model allowable statistical performance specification. 

Prior works in probabilistic Lambert problem \citep{armellin2012high,schumacher2015uncertain,zhang2018covariance,adurthi2020uncertain,gueho2020comparison} have focused on parametric treatments of the problem. These involve approximating either the statistics (e.g., steering moments of endpoint statistics such as covariances) or the dynamics (e.g., Taylor series approximation) or some combinations thereof. We show that significant headway can be made in both theoretical understanding and computational solution of \eqref{ProbLambertFeasibility} and its generalization with stochastic process noise, in the nonparametric setting.

\textbf{Contributions.} By linking problem \eqref{ProbLambertFeasibility} with the Optimal Mass Transport (OMT) theory, we uncover significant insights. \emph{First}, the OMT techniques help proving uniqueness of the solution to \eqref{ProbLambertFeasibility}. \emph{Second}, the methodology used to establish this uniqueness also demonstrates that the identified unique velocity field $\bm{v}$ is optimal w.r.t. certain Lagrangian that involves the gravitational potential $V$. \emph{Third}, this approach facilitates natural extension of \eqref{ProbLambertFeasibility} to scenarios where the velocity includes additive process noise (due to e.g., noisy actuation, stochastic disturbance in atmospheric drag). In such cases, the controlled ordinary differential equation
\begin{align}
\dot{\bm{r}}=\bm{v}(\bm{r},t) 
\label{ODESamplePath}
\end{align}
is substituted with the Itô stochastic differential equation
\begin{align}
\differential\bm{r} = \bm{v}(\bm{r},t)\:\differential t + \sqrt{2\varepsilon}\:\differential\bm{w}(t),
\label{SDESamplePath}
\end{align} 
where $\bm{w}(t)\in\mathbb{R}^3$ denotes the standard Wiener process, and $\varepsilon>0$ represents the strength of the process noise. We show that problem \eqref{ProbLambertFeasibility}, with \eqref{ODESamplePath} replaced by \eqref{SDESamplePath}, evolves into a generalized Schrödinger bridge problem (SBP) \citep{schrodinger1931umkehrung,schrodinger1932theorie,wakolbinger1990schrodinger}. 

\textbf{Organization.} In the remaining of this extended abstract, we summarize our findings as follows. In Sec. \ref{MainResults}, we state our main results with brief discussions clarifying the logical progression of ideas. We eschew the proofs due to space constraints; they can be found in our extended manuscript \citep{teter2024solution}. The closing remarks in Sec. \ref{ConcludingRemarks} include a summary and directions of future work.

\section{Main Results}\label{MainResults}
\textbf{Lambertian Optimal Mass Transport (L-OMT).}\label{sec:LambertOT}
Following \citep{bando2010new,kim2020optimal}, problem \eqref{TPBVP} can be cast as a deterministic optimal control problem:
\begin{subequations}
\begin{align}
&\underset{\bm{v}}{\arg\inf}\quad\displaystyle\int_{t_0}^{t_1}\left(\dfrac{1}{2}\vert \bm{v}\vert ^{2} - V(\bm{r})\right)\differential t \label{LambertOCPobj}\\  
&\dot{\bm{r}} = \bm{v}, \label{LambertOCPdynconstr}\\
& \bm{r}(t=t_0) = \bm{r}_0\:\text{(given)}, \quad \bm{r}(t=t_1) = \bm{r}_{1}\:\text{(given)}\label{LambertOCPterminalconstr}.
\end{align}
\label{LambertOCP}
\end{subequations}
For the probabilistic Lambert problem \eqref{ProbLambertFeasibility}, we need to replace \eqref{LambertOCPterminalconstr} with \eqref{ProbLambertEndpoint}. 

Additionally, the uncertainty in the initial condition $\bm{r}(t=t_0)\sim\rho_0$ changes over time as it follows the path defined by \eqref{LambertOCPdynconstr}. This change in uncertainty is described using the Liouville partial differential equation (PDE)
\begin{align}
\dfrac{\partial\rho}{\partial t} + \nabla_{\bm{r}}\cdot \left(\rho\bm{v}\right)=0, \quad \rho(t=0,\cdot)=\rho_0\:\text{(given)},
\label{LiouvilleIVP}    
\end{align}
where $\rho(\bm{r},t)$ represents the transient joint state PDF, influenced by a feasible control policy $\bm{v}(\bm{r},t)$.
Different control policies $\bm{v}$ in \eqref{LiouvilleIVP} lead to different PDF-valued trajectories, connecting the given $\rho_0,\rho_1\in\mathcal{P}_{2}\left(\mathbb{R}^3\right)$. Therefore, problem \eqref{LambertOCP} with \eqref{LambertOCPterminalconstr} replaced by \eqref{ProbLambertEndpoint}, leads to a generalized OMT formulation
\begin{subequations}
\begin{align}
&\underset{\left(\rho,\bm{v}\right)\in\mathcal{P}_{01}\times\mathcal{V}}{\arg\inf}\quad\displaystyle\int_{t_0}^{t_1}\int_{\mathbb{R}^3}\left(\dfrac{1}{2}\vert \bm{v}\vert ^{2} - V(\bm{r})\right)\rho(\bm{r},t)\:\differential\bm{r}\:\differential t \label{LambertOTobj}\\  
&\dfrac{\partial\rho}{\partial t} + \nabla_{\bm{r}}\cdot \left(\rho\bm{v}\right)=0, \label{LambertOTdynconstr}\\
& \bm{r}(t=t_0) \sim \rho_0\:\text{(given)}, \quad \bm{r}(t=t_1) \sim \rho_{1}\:\text{(given)}\label{LambertOTterminalconstr}.
\end{align}
\label{LambertOT}
\end{subequations}
\noindent where $\mathcal{P}_{01}$ denotes the collection of all PDF-valued trajectories $\rho(\cdot,t)$ that are continuous in $t\in[t_0,t_1]$, supported over $\mathbb{R}^3$, and satisfy $\rho(\cdot,t=t_0)=\rho_0, \rho(\cdot,t=t_1)=\rho_1$. The set $\mathcal{V}$ comprises Markovian finite energy control policies.

We refer to \eqref{LambertOT} as the \emph{Lambertian optimal mass transport (L-OMT)} problem. The classical OMT \citep{benamou2000computational} is its special case $V\equiv 0$. Our first result is that the probabilistic Lambert problem \eqref{ProbLambertFeasibility} admits unique solution that comes with an inverse optimality guarantee.

\begin{theorem}\label{thm:existenceuniqueness} (\textbf{Existence, Uniqueness and Inverse Optimality of Probabilistic Lambert problem}) 
For given $\rho_0,\rho_1\in\mathcal{P}_2\left(\mathbb{R}^{3}\right)$, the the probabilistic Lambert problem \eqref{ProbLambertFeasibility} (equivalently \eqref{LambertOT}) admits a unique solution $\left(\rho^{\rm{opt}},\bm{v}^{\rm{opt}}\right)$ that is, in fact, a minimizer of \eqref{LambertOTobj}.
\end{theorem}
The proof for the above follows from Figalli's theory \citep{figalli2007optimal} for OMTs with cost derived from an action functional.

\textbf{Lambertian Schr\"{o}dinger Bridge (L-SBP).} For a given $\varepsilon>0$ (not necessarily small) and fixed time interval $[t_0,t_1]$ as before, we introduce a stochastic process noise variant of the L-OMT \eqref{LambertOT}, referred hereafter as the \emph{Lambertian Schrödinger Bridge Problem (L-SBP)}, given by
\begin{subequations}
\begin{align}
&\underset{\left(\rho,\bm{v}\right)\in\mathcal{P}_{01}\times\mathcal{V}}{\arg\inf}\quad\displaystyle\int_{t_0}^{t_1}\int_{\mathbb{R}^3}\left(\dfrac{1}{2}\vert \bm{v}\vert^{2} - V(\bm{r})\right)\rho(\bm{r},t)\:\differential\bm{r}\:\differential t \label{LambertSBPobj}\\  
&\dfrac{\partial\rho}{\partial t} + \nabla_{\bm{r}}\cdot \left(\rho\bm{v}\right)=\varepsilon\Delta_{\bm{r}}\rho, \label{LambertSBPdynconstr}\\
& \bm{r}(t=t_0) \sim \rho_0\:\text{(given)}, \quad \bm{r}(t=t_1) \sim \rho_{1}\:\text{(given)}\label{LambertSBPterminalconstr}.
\end{align}
\label{LambertSBP}
\end{subequations}
The dynamic constraint \eqref{LambertSBPdynconstr} replaces \eqref{LambertOTdynconstr}, which corresponds to replacing the sample path dynamics \eqref{ODESamplePath} by \eqref{SDESamplePath}. Our next result is the following.
\begin{theorem}\label{ThmLSBPExistenceUniqueness} 
\textbf{(Uniqueness of L-SBP Solution)} 
Given $\rho_0,\rho_1\in\mathcal{P}_2(\mathbb{R}^3)$ and $\varepsilon>0$, the L-SBP \eqref{LambertSBP}, which is the probabilistic Lambert problem with process noise, admits a unique minimizing pair $(\rho_{\varepsilon}^{\text{opt}},\bm{v}_{\varepsilon}^{\text{opt}})$.
\end{theorem}
The main idea behind the proof for the above is to recast the L-SBP \eqref{LambertSBP} as a relative entropy minimization problem w.r.t. certain Gibbs measure over the path space, i.e., to derive a large deviation principle \citep{dembo2009large} for \eqref{LambertSBP} where the Kullback-Leibler divergence plays the role of the rate functional.

The next result deduces the necessary conditions for optimality for the L-SBP \eqref{LambertSBP}. The corresponding conditions for the L-OMT \eqref{LambertOT} are obtained by specializing the following for $\varepsilon=0$.

\begin{proposition}(\textbf{Conditions for optimality for L-SBP})\label{PropCondOptimalityLSBP}
The optimal solution pair $(\rho_{\varepsilon}^{\text{opt}},\bm{v}_{\varepsilon}^{\text{opt}})$ for the L-SBP \eqref{LambertSBP} solves the system of coupled nonlinear PDEs
\begin{subequations}
    \begin{align}
        &\frac{\partial\psi}{\partial t} + \frac{1}{2}|\nabla_{\bm{r}} \psi|^2 +  \varepsilon\Delta_{\bm{r}}\psi = - V(\bm{r}),\label{HJBSBP}\\
        &\frac{\partial\rho_{\varepsilon}^{\text{opt}}}{\partial t} + \nabla_{\bm{r}}\cdot (\rho_{\varepsilon}^{\text{opt}}\nabla_{\bm{r}} \psi)=\varepsilon\Delta_{\bm{r}}\rho_{\varepsilon}^{\text{opt}}, \label{FPKSBP}
    \end{align}
    \label{FirstOrderConditions}
\end{subequations}
with boundary conditions 
\begin{align}
    \rho_{\varepsilon}^{\rm{opt}}(\bm{r},t=t_0) = \rho_0(\bm{r}), \quad \rho_{\varepsilon}^{\rm{opt}}(\bm{r},t=t_1) = \rho_1(\bm{r}),
    \label{BoundaryConditionsForProp1}
\end{align}
where $\psi\in C^{1,2}\left(\mathbb{R}^{3};[t_0,t_1]\right)$, and the optimal control 
\begin{align}
    \bm{v}_{\varepsilon}^{\rm{opt}} = \nabla_{\bm{r}} \psi(\bm{r}, t).
    \label{inputSBP}
\end{align}
\end{proposition}
In the subsequent Theorem, we employ the Hopf-Cole transformation \citep{hopf1950partial,cole1951quasi} to transcribe \eqref{FirstOrderConditions}-\eqref{BoundaryConditionsForProp1} into a system of \emph{linear} reaction-diffusion PDEs.
\begin{theorem}\label{ThmHopfColeReactionDiffusion}
\textbf{(Linear Reaction-Diffusion PDEs)} For the L-SBP \eqref{LambertSBP} with conditions specified in Theorem \ref{ThmLSBPExistenceUniqueness}, let $(\rho^{\text{opt}}_\varepsilon, \psi)$ be the solution to \eqref{FirstOrderConditions}-\eqref{BoundaryConditionsForProp1}. Define the Schr\"{o}dinger factors $(\widehat{\varphi}_{\varepsilon},\varphi_{\varepsilon})$ via Hopf-Cole transformation $(\rho^{\text{opt}}_\varepsilon, \psi) \mapsto(\widehat{\varphi}_{\varepsilon},\varphi_{\varepsilon})$ as
\begin{subequations}
    \begin{align}
        \varphi_{\varepsilon} &= \exp \left( \frac{\psi}{2 \varepsilon} \right), \label{phi_epsilon}\\
        \widehat{\varphi}_{\varepsilon} &= \rho^{\text{opt}}_\varepsilon \exp \left( - \frac{\psi}{2 \varepsilon} \right).\label{phi_hat_epsilon}
    \end{align} 
\end{subequations}
Then $(\widehat{\varphi}_{\varepsilon},\varphi_{\varepsilon})$ solve linear reaction-diffusion PDEs:
\begin{subequations}
    \begin{align}
        \frac{\partial\widehat{\varphi}_{\varepsilon}}{\partial t} &= \left(\varepsilon\Delta_{\bm{r}} + \frac{1}{2\varepsilon} V(\bm{r})\right)\widehat{\varphi}_{\varepsilon},\label{FactorPDEForward}\\
        \frac{\partial\varphi_{\varepsilon}}{\partial t} &= -\left(\varepsilon\Delta_{\bm{r}} + \frac{1}{2\varepsilon}V(\bm{r})\right)\varphi_{\varepsilon},\label{FactorPDEBackward}
    \end{align}
    \label{FactorPDEs}
\end{subequations}
with coupled boundary conditions
\begin{align}
    &\widehat{\varphi}_{\varepsilon}(\cdot,t=t_0)\varphi_{\varepsilon}(\cdot,t=t_0) = \rho_0, \nonumber\\ &\widehat{\varphi}_{\varepsilon}(\cdot,t=t_1)\varphi_{\varepsilon}(\cdot,t=t_1) = \rho_1.\label{factorBC}    
\end{align}
The optimal pair for \eqref{LambertSBP} is recovered from these factors as
\begin{subequations}
    \begin{align}
        \rho_{\varepsilon}^{\text{opt}}(\cdot,t) &= \widehat{\varphi}_{\varepsilon}(\cdot,t)\varphi_{\varepsilon}(\cdot,t),\label{OptimallyControlledJointPDF}\\
        \bm{v}_{\varepsilon}^{\text{opt}}(t,\cdot) &= 2 \varepsilon\nabla_{(\cdot)}\log\varphi_{\varepsilon}(\cdot,t).\label{OptimalVelocity}
    \end{align}
\end{subequations}
\end{theorem}

\begin{figure}[t] 
\centering
\includegraphics[width=.23\textwidth]{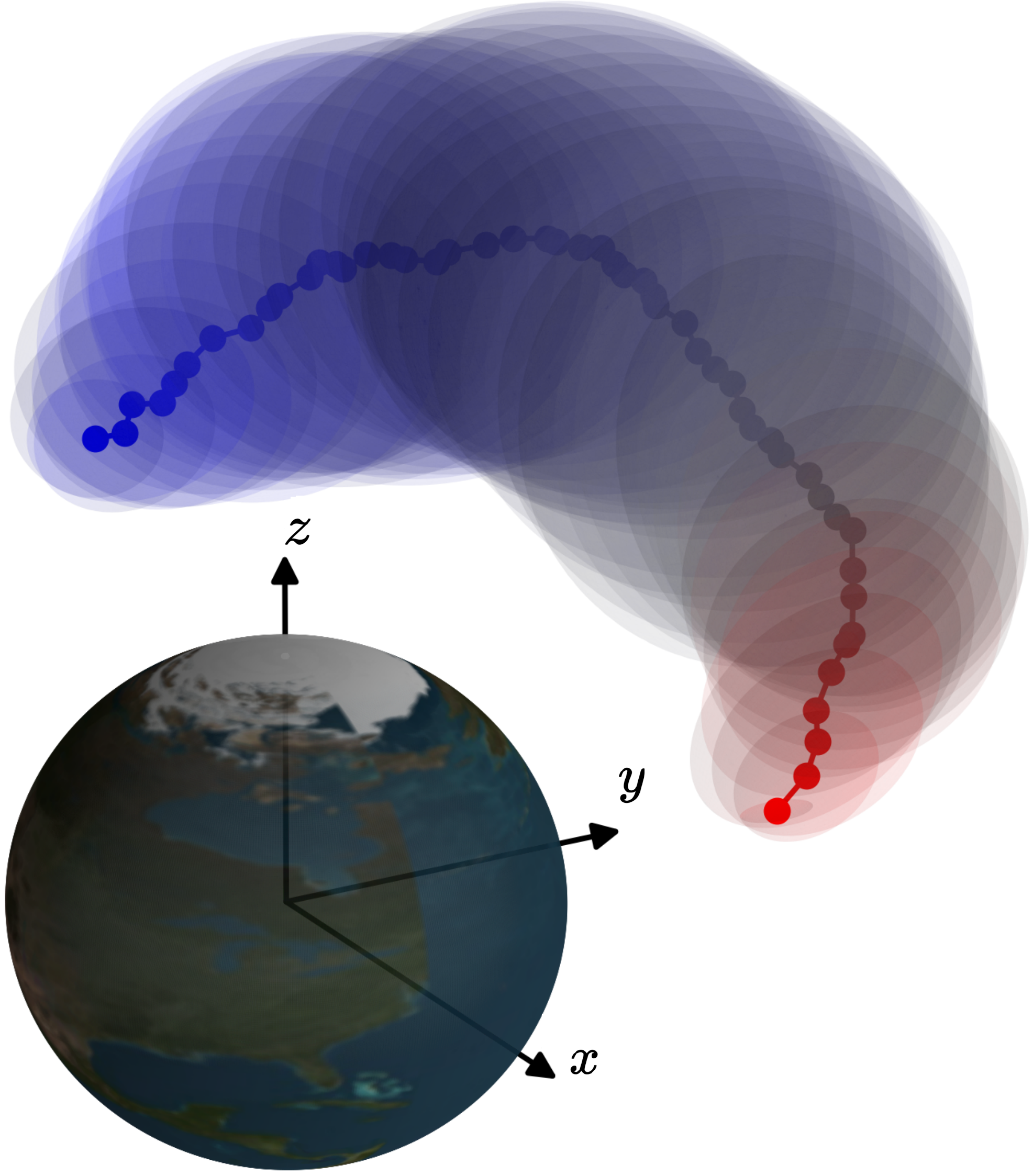}
\caption{Solution of the L-SBP \eqref{LambertSBP}, i.e., probabilistic Lambert problem with process noise, for a numerical case study detailed in \citep[Sec. VI]{teter2024solution}.}
\label{fig:3D_plot}
\end{figure}

\textbf{Computation.}
We propose to solve \eqref{FactorPDEs}-\eqref{factorBC} through a recursive algorithm as follows. Notice that we can solve the initial value problems for PDEs \eqref{FactorPDEs} either via Feynman-Kac path integral or via Fredholm integral \citep[Sec. V]{teter2024solution}. However, the endpoint Schr\"{o}dinger factors $\widehat{\varphi}_{\varepsilon,0}(\cdot) := \widehat{\varphi}_{\varepsilon}(\cdot,t=t_0)$ and $\varphi_{\varepsilon,1}(\cdot) := \varphi_{\varepsilon}(\cdot,t=t_1)$ are not known \emph{a priori}. Our proposed method begins with an everywhere positive guess for $\widehat{\varphi}_{\varepsilon,0}$ and uses it to predict $\widehat{\varphi}_{\varepsilon}(\cdot,t=t_1)$ by integrating \eqref{FactorPDEForward} forward in time. We then estimate $\varphi_{\varepsilon,1}(\cdot)$ at $t=t_1$ using boundary condition \eqref{factorBC}, and integrate \eqref{FactorPDEBackward} backward in time to predict $\varphi_{\varepsilon}(\cdot,t=t_0)$. Then, applying the boundary condition \eqref{factorBC} at $t=t_0$ gives us a new estimate for $\varphi_{\varepsilon,0}(\cdot)$, completing one iteration of the recursion. We repeat this process until it converges. 
In practical numerical simulations, $\rho_0,\rho_1$ are compactly supported, fulfilling the requirement $\rho_0,\rho_1\in\mathcal{P}_2(\mathbb{R}^3)$. The proposed recursion for such compactly supported endpoint data is known \citep[Sec. III]{chen2016entropic} to be contractive w.r.t. Hilbert's projective metric \citep{hilbert1895gerade,bushell1973hilbert,franklin1989scaling}, i.e., enjoys \emph{guaranteed linear convergence}. For further details on this contractive recursion, we refer the readers to \citep{de2021diffusion,pavon2021data,caluya2021reflected,caluya2021wasserstein,10293168}.

Fig. \ref{fig:3D_plot} depicts the solution for an instance of the L-SBP \eqref{LambertSBP} using the computational procedure outlined above. For details, we refer the readers to \citep[Sec. VI]{teter2024solution}.


\section{Concluding Remarks}\label{ConcludingRemarks}
We found that the probabilistic Lambert problem admits unique solution that in fact comes with an optimality certificate. We showed that the same problem with process noise also admits a unique solution that can be computed via certain recursion involving a pair of boundary-coupled linear reaction-diffusion PDEs. Our future work will explore directly computing the Green's functions for these PDEs for improved computation.    

\bibliography{ifacconf}             
          
\end{document}